\documentclass[journal]{IEEEtran}
 \usepackage{graphicx}

\usepackage{amsmath,amsthm}

\usepackage{amsfonts}

\usepackage{cases}
\usepackage{subeqnarray}
\usepackage{subfigure}
\usepackage{url}

\ifCLASSINFOpdf
\else
\fi
\begin{document}
%
\title{\LARGE Charging Scheduling of Electric Vehicles with Local Renewable Energy under Uncertain Electric Vehicle Arrival and Grid Power Price}


\author{Tian~Zhang,
Wei~Chen,~\IEEEmembership{Member,~IEEE,}
        Zhu Han,~\IEEEmembership{Senior Member,~IEEE,}
        and~Zhigang~Cao,~\IEEEmembership{Senior Member,~IEEE}
\thanks{ T. Zhang is with the School of Information Science and Engineering, Shandong University, Jinan 250100, China.
He is also with
Tsinghua University.
E-mail: tianzhang.ee@gmail.com}
\thanks{ W. Chen and Z. Cao are with the State Key Laboratory on
Microwave and Digital Communications, Department of Electronic Engineering,
Tsinghua National Laboratory for Information Science and Technology
(TNList), Tsinghua University, Beijing 100084, China.
E-mail: \{wchen, czg-dee\}@tsinghua.edu.cn}
\thanks{ Z. Han is with the Department of Electrical and Computer Engineering, University
of Houston, Houston, TX 77004, USA.
E-mail: zhan2@uh.edu }
}


%


\maketitle

\begin{abstract}
In the paper, we consider delay-optimal charging scheduling of the electric vehicles (EVs) at a charging station with multiple charge points. The charging station is equipped with renewable energy generation devices and can also buy energy from power grid. The uncertainty of the EV arrival, the
intermittence of the renewable energy, and the variation of the grid power price are taken into account and described as independent Markov processes. Meanwhile, the charging energy for each EV is random. The goal is to minimize the mean waiting time of EVs under the long term constraint on the cost.
We propose queue mapping to convert the EV queue to the charge demand queue and prove the equivalence between the minimization of the two queues' average length. Then we focus on the minimization for the average length of the charge demand queue under long term cost constraint. We propose a framework of Markov decision process (MDP) to investigate this scheduling problem.
The system state includes the charge demand queue length, the charge demand arrival, the energy level in the storage battery of the renewable energy, the renewable energy arrival, and the grid power price. Additionally the number of charging demands and the allocated energy from the storage battery compose the two-dimensional policy.
We derive two necessary conditions of the optimal policy. Moreover, we discuss the reduction of the two-dimensional policy to be the number of charging demands only. We give the sets of system states for which charging no demand and charging as many demands as possible are optimal, respectively.
Finally we investigate the proposed radical policy and conservative policy numerically.
\end{abstract}

\begin{IEEEkeywords}
Electric vehicle, charging scheduling, renewable energy, Markov decision process.
\end{IEEEkeywords}

%
\IEEEpeerreviewmaketitle

\section{Introduction}
As an important method of operation to mitigate the shortage of the fossil fuel and severe environmental problems,
the electric vehicle (EV) technology has attracted much interest in recent years. Compared to conventional vehicles, EVs have advantages in the following aspects: energy efficiency, eco-effect, performance benefits, and energy independence \cite{Infocomworkshop12:Y. Li R. Kaewpuang P. Wang D. Niyato and Z. Han}. However, a fuel driven vehicle can produce less CO2
than an EV if the charging energy is entirely produced by coal-fired power plants \cite{IEVC12:J. Keiser M. Lutzenberger S. Albayrak}. Thus, the renewable energy (e.g., solar or wind energy \cite{ProcIEEE11:L. Xie P. M. S. Carvalho L. A. F. M. Ferreira J. Liu B. Krogh N. Popli and M. D. Ilic}) should be the energy source of the EVs fully or at least partially to achieve the real environmental advantages.
\par
Since EVs are propelled by an electric motor (or motors) that is powered by rechargeable battery packs, EVs need to be charged periodically. Then the EV charging becomes an important topic \cite{Book12:E. Hossain Z. Han and V. Poor,IPEC07:G. B. Shrestha S. G. Ang and S. G. Ang}. In particular, there are some works on the scheduling of EV charging in literature \cite{CCPE10:O. Sundstrom and C. Binding}-\cite{TPS12:C.-T. Li C. Ahn H. Peng and J. S. Sun}.
\par
In \cite{CCPE10:O. Sundstrom and C. Binding}, the EV battery charging behavior was optimized with the objective to minimize charging costs, achieve satisfactory state-of-energy levels, and
optimal power balancing. In \cite{ACC10:S. Bashash S. J. Moura and H. K. Fathy}, the problem of optimizing plug-in hybrid electric vehicle (PHEV) charge trajectory (i.e. timing and rate of the charging) was studied to reduce the energy cost and battery degradation. For the purpose of improving the satisfiability of EVs, a reservation-based scheduling algorithm for the charging station to decide the service order of multiple requests was proposed in \cite{SuComs10:H.-J. Kim J.Lee G.-L. Park M.-J. Kang and M. Kang}.
In \cite{PESGM11:S. Sojoudi and S. H. Low}, a joint optimal power flow (OPF)-charging (dynamic) optimization problem was formulated with the goal of minimizing the generation and charging costs while satisfying the network, physical and inelastic-load constraints. In \cite{PESGM11:W. Su and M.Y. Chow}, utilizing the particle swarm
optimization, a proposed algorithm optimally manages a large number of PHEVs charging at a municipal parking station.
In \cite{VANET11:H. Qin and W. Zhang}, the minimization of the waiting time for EV
charging via scheduling charging activities spatially and temporally in a large-scale road network was investigated.
By modeling an EV charging system as a cyber-physical system, a decentralised online EV charging scheduling scheme was developed in \cite{IJPEDS12:R. Jin B. Wang P. Zhang and P.B. Luh}. In \cite{TPS12:L. Gan U. Topcu and S. H. Low}, the authors formulated the EV
charging scheduling problem to fill the electric load valley as an optimal control problem, and a decentralized algorithm was derived. In \cite{TCST12: Z. Ma D. S. Callaway and I. A. Hiskens}, a strategy to coordinate the
charging of plug-in EVs (PEVs) was proposed by using the non-cooperative games \cite{MSP12:W. Saad Z. Han H. V. Poor and T. Basar}. Flexible charging optimization for EVs considering distribution grid constraints, both voltage and power, was investigated in \cite{TSG12:O. Sundstrom and C. Binding}.
In \cite{ACC12:J. Huang V. Gupta and Y.-F. Huang}, the trade off between distribution system load with
quality of charging service was considered, and the centralized
algorithms to schedule the charging of vehicles were designed.
In \cite{ACC12:A. Subramanian M. Garcia A. Dominguez-Garcia D. Callaway K. Poolla and P. Varaiya} and \cite{SmartGridComm12:S. Chen and L. Tong}, real-time scheduling policies of EV charging were considered when both the renewable energy and energy
from the grid are available. In \cite{TPS12:C.-T. Li C. Ahn H. Peng and J. S. Sun}, the PEV charging and wind power scheduling were integrated, and the synergistic control algorithm of plug-on vehicle charging and wind power scheduling was proposed.
\par
In the paper, we focus on the scheduling approach of EV charging at a charging station.
The charging station has multiple charge points and is equipped with renewable energy generation devices and storage battery.
The charged energy at a charge point during a period is constant and is called an energy block.
We model the arrival of the renewable energy as a Markov chain. The charging energy can also be purchased from power grid, and the price changes also according to another Markov chain. The arrival of the EVs is assumed as a Markov process. Once an EV arrives at the charging station, it waits in a queue before charging. In each period, the charging station chooses some EVs from the head of the queue for charging. Meanwhile, the station also determines
how much energy is supplied from the storage battery (the rest of the required energy is supplied from the power grid). The objective is minimizing the mean waiting time of EVs under the long term cost constraint.
\par
Since the amount of charging energy (i.e., the number of energy blocks to charge) for each EV is random, the scheduling problem is very challenging. We propose queue mapping method to solve the difficulty. We map the EV queue to a charge demand queue. In the charge demand queue, each demand means an energy block that need to charge and some consecutive demands correspond to an EV's required charing energy. We prove that the minimization of the average EV queue length is equivalent to the minimization of the average charge demand queue length. Then we focus on the charge demand queue minimization under the cost constraint.
The scheduling problem can be equivalently reconstructed as follows.
The demand arrives according to a discrete-time batch Markovian arrival process (D-BMAP) and waits in the charge demand queue before service (charging). In each period, the charging station chooses some demands from the head of the charge demand queue for charging. Meanwhile, the station also determines
how much energy is supplied from the storage battery (the rest of the required energy is supplied from the power grid). The objective is minimizing the mean length of the charge demand queue under the long term cost constraint.
\par
Next, we find that the reconstructed optimization problem can be studied under a Markov decision process (MDP) framework.
The system state contains the charge demand queue length, the demand arrival, the energy level in the storage battery of the renewable energy, the renewable energy arrival, and the grid power price. Meanwhile, the number of charging demands and the allocated energy from the storage battery constitute the two-dimensional policy. We find that the general case of the reconstructed optimization problem can be analyzed similarly as the analysis of a special case. Then we focus on the analysis of the special case that is formulated as a constrained MDP \cite{Book99:E. Altman}.
We analyze the optimal two-dimensional policy of the constrained MDP by transforming to an average cost MDP and its corresponding discount cost MDP thereafter. First, the constrained MDP is converted to an unconstrained MDP by using Lagrangian relaxation. Moreover, we derive that the optimal solution of the unconstrained MDP with a certain Lagrangian multiplier is the optimal for the original constrained MDP. Next, the unconstrained MDP can be analyzed by transforming to its corresponding discount cost MDP. We obtain two necessary conditions for the optimal solution. Third, we analyze the relations between the two elements of the two-dimensional policy, and find that the number of charging demands\footnote{ In the special case, we can use \lq\lq EV\rq\rq ~and \lq\lq demand\rq\rq ~interchangeably.} is dominant. Thus, we propose a conjecture that the constrained MDP problem can be reduced to a MDP problem with the policy to be the number of charging demands only. We then derive the conditions of the system state when the policy that charging no demand is optimal. We also obtain the system state conditions when charging as many demands as possible is optimal.
\par
The rest of the paper is structured as follows.
In Section II, the system model is described and we formulate an optimization problem that can be studied under the framework of MDP. Section III presents a spacial case of the formulated optimization problem as a constrained MDP to demonstrate the solving process of the general case.
Next, we analyze the optimal policy of the constrained MDP in Section IV. In Section V, the numerical results are performed. Finally, Section VI concludes the paper.

\section{System model and problem formulation}
\theoremstyle{definition} \newtheorem{lemma}{Lemma}
\theoremstyle{definition} \newtheorem{property}{Property}

\begin{figure}[!t]
\centering
\includegraphics[width=3.5in]{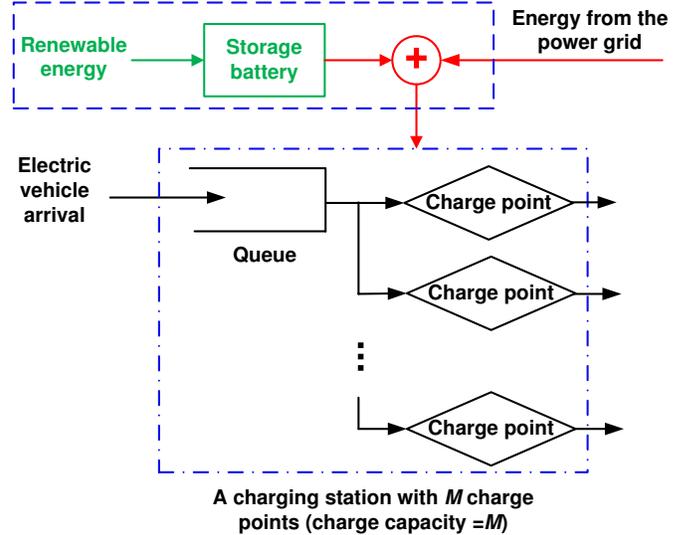}
\caption{System model}
\label{Sysmodel}
\end{figure}

\theoremstyle{definition} \newtheorem{proposition}{Proposition}
\theoremstyle{definition} \newtheorem{conjecture}{Conjecture}
Time is divided into periods of length $\tau$ each. The EVs arrive at the charging station according to a finite-state ergodic Markov chain $\{A[n]\}$.
The EVs wait in a queue before charging as illustrated in Fig. \ref{Sysmodel}. The charging station has $M$ charge points, i.e., at most $M$ EVs can be charged in each period. The charging station has renewable energy generation devices, and it can also gets power from the power grid. The renewable energy is modeled as another finite-state ergodic Markov process $\{E_a[n]\}$.
The renewable energy is viewed as free, and the price for the grid power during the $n$-th period is denoted as $P[n]$.
The grid power price remains static during each period and changes between different periods. The sequence of the price, $\{P[n]\}$, is a finite-state ergodic Markov chain.
We assume that the charged energy at one charging point during a period is constant, and is denoted as $\mathcal{E}$.\footnote{It is assumed that if an EV utilizes $m$ charge points during a period, the amount of charged energy is $m\mathcal{E}$.}
In the $n$-th period, $k[n]$ EVs from the head of the EV queue are allowed to charge.
During the $n$-th period, the charging station allocates $w[n]$ power from the storage battery, and the rest power
will be supplied by the power grid.
Assume that the required charging energy of the EV, $E_c$, is independent on each other, and $E_{c}=L\mathcal{E}$ with $L$ being uniformly distributed in $[1,2,...,C]$,\footnote{$C$ is a given constant.} i.e., $L \sim \mathcal{U}[1,...,C]$.
\par
\begin{figure}[!t]
\centering
\includegraphics[width=3.5in]{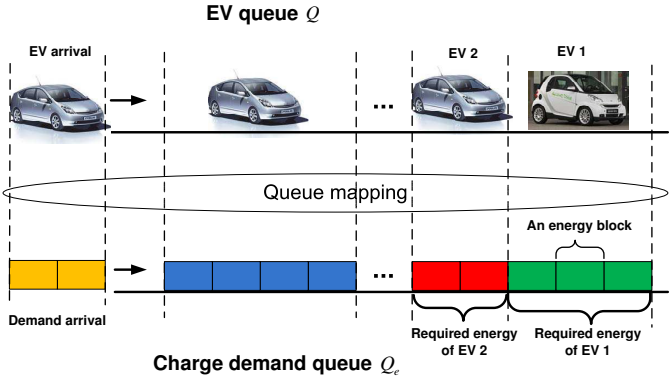}
\caption{Queue mapping}
\label{queue mapping}
\end{figure}

Direct analysis of the EV queue length under the long term cost constraint is difficult due to the randomness of $L$.
We propose the queue mapping method as shown in Fig. \ref{queue mapping}.
Each EV in the EV queue corresponds to several consecutive charge demands (the number of the demands denotes the amount of required energy) in the charge demand queue.\footnote{A demand means $\mathcal{E}$ energy (i.e., an energy block) need to be charged. In Fig. \ref{queue mapping}, the first EV (EV 1) in the EV queue wants to charge $3\times\mathcal{E}$, then it corresponds to the first three consecutive charge demands in the charge demand queue. The second EV (EV 2) charges $2\times\mathcal{E}$, then it corresponds to the two consecutive charge demands after the first EV's corresponding charge demands.} The number of EVs at the beginning of the $n$-th period is $Q[n]$ and the length of the charge demand queue is denoted as $Q_e[n]$.
We convert the average EV queue length minimization to the average charge demand queue minimization. Furthermore, we will prove that they are equivalent.
\par
The demand arrival can be given by
$
A^{'}[n]=\sum_{i=1}^{A[n]}L_i,
$
where $L_i \sim \mathcal{U}[1,...,C]$.
\par
\emph{Remark: As $\{A[n]\}$ is a Markov chain, we can derive that $\{A^{'}[n]\}$ is a D-BMAP.}
\par
In the $n$-th period, $k^{'}[n]$ demands from the head of the charge demand queue are allowed to charge.
During the $n$-th period, the charging station
allocates $w[n]$ power from the storage battery, and the rest power
will be supplied by the power grid.
Denote the number of charged demands in the $n$-th period as
$K^{'}[n]$.
The evolution of the charge demand queue length, $Q_e[n]$, is
$
Q_e[n+1]=Q_e[n]-K^{'}[n]+A^{'}[n].
$
Denote the capacity of the renewable energy storage battery as $E_{max}$. The stored battery energy at the beginning of the $n$-th period is $E_b[n]$.
The battery energy evolution can be expressed as
\begin{eqnarray}\label{battery energy evolution}
E_b[n+1]=\min\big\{E_b[n]-W[n]\tau+E_a[n],E_{max}\big\} \nonumber\\
:=\big(E_b[n]-W[n]\tau+E_a[n]\big)^-.
\end{eqnarray}
The cost in the  $n$-th period is
$
 \mathcal{C}^{'}[n]=\Big(\frac{K^{'}[n]\mathcal{E}}{\tau}-W[n]\Big)^{+}P[n].
 $
\par
Denote the state space as $\mathcal{X}^{'}$ and denote the action space as $\mathcal{A}^{'}$.
Let the (random) system state and action in the $n$-th period be $X^{'}[n]=(Q_e[n],A^{'}[n],E_b[n],E_{a}[n],P[n]) \in \mathcal{X}^{'}$ and $\big(K^{'}[n],W[n]\big) \in \mathcal{A}^{'}$, respectively.
Define a policy $\pi^{'}=(\pi_0^{'},\pi_1^{'},\cdots)$ with $\pi_n^{'}$ generating an action $(k^{'}[n],w[n])$ with a probability \cite{Book99:E. Altman,Book02:E. Feinberg and A. Shwartz} at the $n$-th period. We denote the set of all policies as $\Pi^{'}$. Let $x^{'}[n]=(q_e[n],a^{'}[n],e_b[n],e_a[n],p[n])$ be a (fixed) system state.
The feasible $(k^{'}[n],w[n])$ in state $x^{'}[n]$ belongs to $\mathcal{K^{'}}(x^{'}[n])=\big\{0,1,\cdots,\min\{q_e[n],M\}\big\}\times \mathcal{W}(x^{'}[n])=\{0,\frac{1}{\tau},\cdots,\frac{e_b[n]}{\tau}\}$.\footnote{The energy has been discretized.}
The optimization problem that minimizes the mean charge demand queue length under the long term cost constraint, $\mathop B\limits^\_$, can be expressed as
 \begin{eqnarray} \label{Extended optimization problem}
\min_{\pi \in \Pi^{'}} D_x^{\pi^{'}}  : = \mathop {\lim \sup }\limits_{n \to \infty} \frac{1}{n}\mathbb{E}_{x^{'}}^{\pi^{'}}\left[\sum\limits_{i =0}^{n-1} Q_e[i]\right]
\end{eqnarray}
\begin{subequations}

\begin{numcases}{ \mbox{s.t.}}
B_{x^{'}}^{\pi^{'}}  : = \mathop {\lim \sup }\limits_{n \to \infty} \frac{1}{n}\mathbb{E}_{x^{'}}^{\pi^{'}}\left[\sum\limits_{i=0}^{n-1} {\mathcal{C}^{'}[i]}\right]
 \le \mathop B\limits^\_,\\
K^{'}[i]\le \min\{Q_e[i],M\},\label{original Cons1}\\
W[i]\le \frac{E_{b}[i]}{\tau },\label{original Cons2}
\end{numcases}
\end{subequations}
with initial state $x^{'}=(q_e,a^{'},e_b,e_a,p)$.
\par
Since D-BMAP can be represented by a two-dimensional discrete-time Markov chain (DTMC) \cite{BMAP10:J. D. Cordeiro and J. P. Kharoufeh}, the optimization problem in (\ref{Extended optimization problem}) can be analyzed in the framework of MDP.
Moreover, the following lemma proves the equivalence of the mean energy demand queue length minimization and the mean EV queue length minimization.
\begin{lemma}\label{equivalence energy demand queue and EV queue}
The minimization of the mean charge demand queue length is equivalent to the minimization of the mean EV queue length.
\end{lemma}
\begin{IEEEproof}
See Appendix \ref{proof of equivalence energy demand queue and EV queue}.
\end{IEEEproof}
\section{Simplified problem}
For conciseness, we give a special case of problem (\ref{Extended optimization problem}) in the section and investigate this relatively simplified problem in the following of the paper to show the solving process.
General cases can be analyzed through similar solving process.
\par
When $C=1$, we have $L_i=1$. Then the queue mapping is an identity transform and \lq\lq EV\rq\rq ~and \lq\lq demand\rq\rq ~are interchangeable. Thus, we can directly analyze the EV queue using the MDP framework. We have $K[n]=K^{'}[n]$, $A[n]=A^{'}[n]$ and $Q[n]=Q_e[n]$.
The queue length evolution is
\begin{eqnarray}\label{queue length evolution}
Q[n+1]=Q[n]-K[n]+A[n].
\end{eqnarray}
The battery energy evolution is the same as (\ref{battery energy evolution}).
The cost at the $n$-th period is given by
\begin{eqnarray}\label{stage cost}
\mathcal{C}[n]=\Big(\frac{K[n]\mathcal{E}}{\tau}-W[n]\Big)^{+}P[n],
\end{eqnarray}
where $(\cdot)^+:=\max\{\cdot,0\}$.
The system state becomes $X[n]=\left(Q[n],A[n],E_{b}[n],E_a[n],P[n]\right)$ with state space $\mathcal{X}$
and the action is $\left(K[n],W[n]\right)$ with action space $\mathcal{A}$. $\left\{X[n],\left(K[n],W[n]\right)\right\}$ is a controlled Markov process.
Define a policy $\pi=(\pi_0,\pi_1,\cdots)$ that $\pi_n$ generates an action $(k[n],w[n])$ with a probability at the beginning of the $n$-th period. We denote the set of all policies as $\Pi$.
The feasible $(k[n],w[n])$ in state $x[n]$ belongs to $\mathcal{K}(x[n])=\big\{0,1,\cdots,\min\{q[n],M\}\big\}\times \mathcal{W}(x[n])=\{0,\frac{1}{\tau},\cdots,\frac{e_b[n]}{\tau}\}$.
A stationary deterministic policy is $\pi=(g,g,\cdots)$, where $g$ is a measurable mapping from $\mathcal{X}$ to $\mathcal{K}(x[n]) \times \mathcal{W}(x[n])$.
Our objective is to find a policy that minimizes the mean queue delay under the long run constraint on the cost.
The optimization problem (i.e., the constrained MDP) is given by
\begin{eqnarray} \label{original optimization problem}
\min_{\pi \in \Pi} D_x^\pi  : = \mathop {\lim \sup }\limits_{n \to \infty} \frac{1}{n}\mathbb{E}_{x}^{\pi}\left[\sum\limits_{i =0}^{n-1} Q[i]\right]
\end{eqnarray}
\begin{subequations}
\begin{numcases}{ \mbox{s.t.}}
B_x^\pi  : = \mathop {\lim \sup }\limits_{n \to \infty} \frac{1}{n}\mathbb{E}_{x}^{\pi}\left[\sum\limits_{i=0}^{n-1} {\mathcal{C}[i]}\right]
 \le \mathop B\limits^\_,\\
K[i]\le \min\{Q[i],M\},\label{original Cons1}\\
W[i]\le \frac{E_{b}[i]}{\tau },\label{original Cons2}
\end{numcases}
\end{subequations}
where $x=(q,a,e_b,e_a,p) \in \mathcal{X}$ is the initial system state.
\par
\emph{Remark: (\ref{original optimization problem}) is the special case of (\ref{Extended optimization problem}) with $C=1$. $C=1$ means that EVs charge the same amounts of energy, $\mathcal{E}$
(e.g., an EV production company).
}
\section{Analysis of the optimal policy}
In this section, we perform theoretical study on the optimal policy.
First, we prove that the constrained MDP can be analyzed through an unconstrained MDP. Then, we focus on the analysis of the unconstrained MDP. We analyze the unconstrained MDP by using its corresponding discount MDP. Next, we consider the dimension reduction of the two-dimensional policy. Finally, we propose two stationary deterministic policies based on the theoretical results.
\subsection{Transformation to the unconstrained MDP and discount MDP}
Define $f_{\beta}(x,k,w):=\beta\big(\frac{k\mathcal{E}}{\tau}-w\big)^+p+ q$. We have the following unconstrained MDP (i.e., UP$_\beta$).
\begin{eqnarray}\label{UP_beta}
\min_{\pi} J_\beta ^\pi(x)  : = \mathop {\lim \sup }\limits_{n \to \infty} \frac{1}{n}\mathbb{E}_x^\pi  \left[\sum\limits_{i = 0}^{n - 1} {f_\beta (X[i],K[i],W[i])}\right].
\end{eqnarray}
\par
\emph{Remark: UP$_\beta$ is an average cost MDP. Its optimal solution is referred to as the average cost optimal policy.}
\par
The following lemma reveals that the constrained problem has the same solution as UP$_\beta$ with a certain $\beta$.
\begin{lemma}\label{ConMDP to UnconMDP}
There exists $\beta >0$ for which the optimal solution of the unconstrained MDP in (\ref{UP_beta}) (i.e., UP$_\beta$) is also optimal for the constrained MDP in (\ref{original optimization problem}).
\end{lemma}
\begin{IEEEproof}
See Appendix \ref{proof of conMDP to unconMDP}.
\end{IEEEproof}
Next, we define a discount cost MDP with discount factor $\alpha$ corresponding to UP$_\beta$ for each initial system state $x=(q,a,e_b,e_a,p)$, with value function
\begin{eqnarray}
V_\alpha(x)=\min_{\pi}  \mathbb{E}_x^\pi \left[\sum_{i=0}^{\infty} \alpha^{i}f_\beta \left(X[i],K[i],W[i]\right)\right].
\end{eqnarray}
The optimal solution for the discounted problem is
called a discount optimal policy.
\par
The following lemma reveals the existence of the optimal stationary deterministic  policy of UP$_\beta$, and furthermore, how to derive the average cost optimal policy.
\begin{lemma}\label{discout to UP}
There exists a stationary deterministic policy $(k,w)$ that solves UP$_\beta$, which can be obtained as a limit of
discount optimal policies as $\alpha \to 1$.
\end{lemma}
\begin{IEEEproof}
See Appendix \ref{proof of discout to UP}.
\end{IEEEproof}
\par
Based on the above analysis, we find that the constrained MDP can be analyzed through the defined average cost MDP and its corresponding discount cost MDP thereafter.
Hence, we first investigate the solution of the discount cost MDP in the following subsection.
\subsection{The discount optimal policy}
For state-action pair $\left(x=(q,a,e_b,e_a,p),(k,w)\right)$,
let $u=q-k$ and $\eta=e_b-w\tau$. Then $(u(x),\eta(x))$ can also define a stationary deterministic policy.
Then, the discounted cost optimality equation \cite{Book96:O. H. Lerma and J. B. Lassere,Book02:E. Feinberg and A. Shwartz} is given by
\begin{eqnarray}\label{DCOE}
 \lefteqn{
V_{\alpha}(q,a,e_b,e_a,p)=\min_{
\footnotesize
\begin{array}{c}
u\in \{0,1,\cdots,\min\{q,M\}\} \\
  \eta \in \{0,1,\cdots,e_b\}
\end{array}
}
}
\nonumber\\
&&\bigg\{\beta  \Big(\frac{(q-u)\mathcal{E}}{\tau}-\frac{e_b-\eta}{\tau}\Big)^+p+q  \nonumber\\
&+& \alpha \mathbb{E}_{a,e_a,p}\left[V_{\alpha}(u+A,A,(\eta+E_a)^-,E_a,P)\right]
   \bigg\},
\end{eqnarray}
and the corresponding value iteration algorithm (or successive approximation method) is
\begin{eqnarray}\label{iteration for V_alpha}
 \lefteqn{
V_{\alpha,n}(q,a,e_b,e_a,p)=\min_{
\footnotesize
\begin{array}{c}
  u\in \big\{0,1,\cdots,\min\{q,M\}\big\} \\
  \eta \in \{0,1,\cdots,e_b\}
\end{array}}
}
\nonumber\\
&&\bigg\{ \beta \Big(\frac{(q-u)\mathcal{E}}{\tau}-\frac{e_b-\eta}{\tau}\Big)^+p+q +  \alpha\times \nonumber\\
&&\mathbb{E}_{a,e_a,p}\left[V_{\alpha,n-1}(u+A,A,(\eta+E_a)^-,E_a,P)\right]
   \bigg\}
\end{eqnarray}
 with $V_{\alpha,0}(q,a,e_b,e_a,p)=0.$
 \par
First, regarding $V_{\alpha}(q,h,a,e_b,e)$, we have the following properties (Property \ref{Increasing in q} - Property \ref{convexity}).
\begin{property}\label{Increasing in q}
$V_{\alpha}(q,h,a,e_b,e)$ is an increasing
function of $q$.
\end{property}
\begin{IEEEproof}
See Appendix \ref{proof of Increasing in q}.
\end{IEEEproof}
\begin{property}\label{Non-increasing in eb}
$V_{\alpha}(q,a,e_b,e_a,p)$ is a non-increasing function of $e_b$.
\end{property}
\begin{IEEEproof}
See Appendix \ref{proof of Non-increasing in eb}.
\end{IEEEproof}
In practice, the allocated renewable energy will not surpass the required charging energy. Thus, $k\mathcal{E} \ge w\tau$, i.e.,
\begin{eqnarray}\label{Assumption1}
\frac{(q-u)\mathcal{E}}{\tau}-\frac{e_b-\eta}{\tau}\ge 0.
\end{eqnarray}
\begin{property}\label{convexity}
$V_{\alpha}(q,a,e_b,e_a,p)$ is convex in $(q,e_b)$.
\end{property}
\begin{IEEEproof}
See Appendix \ref{proof of convexity}.
\end{IEEEproof}
Next, the following two lemmas reveal two necessary conditions for the optimality, respectively.
\begin{lemma}\label{Sufficient condition for the non-optimality}
 In state $x=(q,a,e_b,e_a,p)$, $(u(x),\eta(x))$ is not the discount optimal solution if $u(x)> q-\min\{q,M\}$ and $\eta(x)+e_a>E_{max}$.
\end{lemma}
\emph{Remark: Lemma \ref{Sufficient condition for the non-optimality} reveals the sufficient condition for the non-optimality, and it can also be viewed as the necessary condition for the optimality. That is to say, any optimal solutions should not satisfy the condition. }
\begin{lemma}\label{Necessary condition for the discout optimality}
Denote the discount optimal policy in state $x=(q,a,e_b,e_a,p)$ as $(u^*(x),\eta^*(x))$. Then, $(u^*(x),\eta^*(x))$
satisfies the following inequality array\footnote{Using Property \ref{convexity}, we can derive that
$Z_1(u,a,\eta,e_a,p)\le Z_1(u+1,a,\eta,e_a,p)$,
$Z_2(u,a,\eta,e_a,p)\le Z_2(u,a,\eta+1,e_a,p)$, and  $Z_3(u,a,\eta,e_a,p)\le Z_3(u+1,a,\eta+1,e_a,p)$.}
\begin{eqnarray}\label{u}
Z_1(u^*,a,\eta^*,e_a,p) \le \beta \frac{\mathcal{E}}{\tau}p
 \le Z_1(u^*+1,a,\eta^*,e_a,p),
\end{eqnarray}
\begin{eqnarray}\label{eta}
Z_2(u^*,a,\eta^*,e_a,p)\le \beta \frac{-p}{\tau}
\le Z_2(u,a,\eta^*+1,e_a,p),
\end{eqnarray}
\begin{eqnarray}\label{u eta}
\lefteqn{
Z_3(u^*,a,\eta^*,e_a,p) \le \beta \frac{p}{\tau}(\mathcal{E}-1)
}
 \nonumber\\
 &\le& Z_3(u^*+1,a,\eta^*+1,e_a,p),
\end{eqnarray}
where
\begin{eqnarray}
\lefteqn{
Z_1(u,a,\eta,e_a,p)
}\nonumber\\
&=&\alpha \mathbb{E}_{a,e_a,p}\big[G_{1}(u+A,A,(\eta+E_a)^-,E_a,P)\big]
\end{eqnarray}
with
\begin{eqnarray}
  \lefteqn{
G_1(q,a,e_b,e_a,p)=
}
\nonumber\\
&&
V_{\alpha}(q,a,e_b,e_a,p)-V_{\alpha}(q-1,a,e_b,e_a,p),
\end{eqnarray}
\begin{eqnarray}
\lefteqn{
Z_2(u,a,\eta,e_a,p)
}
\nonumber\\
&=&\alpha \mathbb{E}_{a,e_a,p}\big[V_\alpha(u+A,A,(\eta+E_a)^-,E_a,P)
\nonumber\\
&-&V_\alpha(u+A,A,(\eta-1+E_a)^-,E_a,P)\big],
\end{eqnarray}
and
\begin{eqnarray}
\lefteqn{
Z_3(u,a,\eta,e_a,p)=
}
\nonumber\\
&&
\alpha \mathbb{E}_{a,e_a,p}\big[V_{\alpha}(u+A,A,(\eta+E_a)^-,E_a,P)
\nonumber\\
&-&V_{\alpha}(u-1+A,A,(\eta-1+E_a)^-,E_a,P)\big].
\end{eqnarray}
\end{lemma}
\begin{IEEEproof}
See Appendix \ref{proof of Necessary for the discout}.
\end{IEEEproof}
\emph{Remark: Lemma \ref{Necessary condition for the discout optimality} gives the necessary condition of the discount optimality, i.e., the optimal policy (or policies) should be the solution(s) of the inequality array. Specially,
if the inequality array has a single solution, the corresponding single solution is the optimal policy since the existence of the optimal policy. }
\subsection{The average cost optimal policy}
First, Lemma \ref{Sufficient condition for the non-optimality} still holds for the average cost MDP.
Next, based on Lemma \ref{discout to UP} and Lemma \ref{Necessary condition for the discout optimality}, we have the following lemma.
\begin{lemma}\label{Necessary condition for the beta optimality}
 Given state $x=(q,a,e_b,e_a,p)$, the average cost optimal policy $(u^*(x),\eta^*(x))$
 should satisfy the following inequality array
\begin{eqnarray}\label{u beta}
\tilde{Z}_1(u^*,a,\eta^*,e_a,p) \le \beta \frac{\mathcal{E}}{\tau}p
 \le \tilde{Z}_1(u^*+1,a,\eta^*,e_a,p),
\end{eqnarray}
\begin{eqnarray}\label{eta beta}
\tilde{Z}_2(u^*,a,\eta^*,e_a,p)\le \beta \frac{-p}{\tau}
\le \tilde{Z}_2(u,a,\eta^*+1,e_a,p),
\end{eqnarray}
\begin{eqnarray}\label{u eta beta}
\lefteqn{
\tilde{Z}_3(u^*,a,\eta^*,e_a,p) \le \beta \frac{p}{\tau}(\mathcal{E}-1)
}
 \nonumber\\
 &\le& \tilde{Z}_3(u^*+1,a,\eta^*+1,e_a,p),
\end{eqnarray}
where
$$\tilde{Z}_1(u,a,\eta,e_a,p)=\lim\limits_{\alpha \to 1}Z_1(u,a,\eta,e_a,p),$$ $$\tilde{Z}_2(u,a,\eta,e_a,p)=\lim\limits_{\alpha \to 1}Z_2(u,a,\eta,e_a,p),$$ and $$\tilde{Z}_3(u,a,\eta,e_a,p)=\lim\limits_{\alpha \to 1}Z_3(u,a,\eta,e_a,p).$$
\end{lemma}
\subsection{Reducing the policy's dimension}
The number of charging EVs $k$ and the power allocation from the battery $w$ are coupled together, they affect each other.
However, if we assume that $k$ has been chosen, then the required total power has been fixed. In this case, we will allocate as much power as possible from the battery to meet the required total power, i.e., the greedy policy for the battery power allocation. This is because the power from the battery is free (please refer to (\ref{stage cost})). We can guess that the greedy allocation strategy of battery power is the optimal policy. However, it is difficult to prove. The difficulty lies in the fact that the remaining battery energy will affect the future action and cost (e.g., (\ref{DCOE})).
On the other hand, once $w$ has been fixed, the power allocation from the power grid can also affect $k$. In summary, when $k$ is chosen, the optimal $w^*$ is the greedy policy. By contrast, if $w$ is fixed, the optimal $k$ is not fixed, we need to solve the power allocation from the power grid to find the optimal $k^*$. Thus, we can reduce the policy from $(k,w)$ to $k$. We have the following conjecture.
\begin{conjecture}\label{dimensionreduction}
Let $\pi_k=(k[0],k[1],\cdots)$, and (\ref{original optimization problem}) can be converted to
\begin{eqnarray} \label{original optimization problem converted to policy only r}
\min_{\pi_k} B_x^{\pi_k}  : = \mathop {\lim \sup }\limits_{n \to \infty}\frac{1}{n}\mathbb{E}_{x}^{\pi_k}\left[\sum\limits_{i =0}^{n-1} {Q[i]}\right]
\end{eqnarray}
\begin{subequations}
\begin{numcases}{\mbox{s.t.}}
B_x^{\pi_k}  : = \mathop {\lim \sup }\limits_{n \to \infty} \frac{1}{n}\mathbb{E}_{x}^{\pi_k}\Big[\sum\limits_{i=0}^{n-1} \Big(\frac{K[i]\mathcal{E}}{\tau}\nonumber\\
\qquad -\frac{1}{\tau}\min\left\{K[i]\mathcal{E},E_{b}[i]\right\}\Big)^{+}P[i]\Big]
  \le \mathop B\limits^\_,\\
 K[i] \le \min\{Q[i],M\},
\end{numcases}
\end{subequations}
where the evolution of energy in the battery becomes
\begin{eqnarray}
E_{b}[i+1]=
\left(E_b[i]-\min\left\{K[i]\mathcal{E},E_{b}[i]\right\}+E_a[i]\right)^{-}.
\end{eqnarray}
\end{conjecture}
\emph{Remark:
The policy can be reduced in dimension ($(k,w) \to k$).
If the stated $\beta$ in Lemma \ref{ConMDP to UnconMDP} satisfying $\beta \gg 1$, Conjecture \ref{dimensionreduction} can be proved based on (\ref{DCOE}) in addition with Lemma \ref{discout to UP} and Lemma \ref{ConMDP to UnconMDP}.}
\par
In the following, we discuss the optimal policy after dimension reduction.
For state-action pair $\left(x=(q,a,e_b,e_a,p),k\right)$,
let $u=q-k$, and $u(x)$ can also define a stationary deterministic policy.
 We have the following lemmas to reveal the properties of the optimal policy.
\begin{lemma}\label{Necessary condition for the discout optimality of u}
Denote the discount optimal policy in state $x=(q,a,e_b,e_a,p)$ as $u^*(x)$. Then, $u^*(x)$
satisfies
\begin{eqnarray}
Z(u^*)\le \beta\frac {\mathcal{E}}{\tau}p\le Z(u^*+1),
\end{eqnarray}
where
\begin{eqnarray}
\lefteqn{
Z(u)=\alpha \mathbb{E}_{a,e_a,p}\Big[V_{\alpha}(u+A,A,(\eta(u)+E_a)^-,E_a,P)
}  \nonumber\\
&-&V_{\alpha}(u-1+A,A,(\eta(u-1)+E_a)^-,E_a,P)\Big]
\nonumber\\
&+& \beta\frac {\eta(u)-\eta(u-1)}{\tau}p
\end{eqnarray}
with
\begin{eqnarray}
 \eta(u):=\max\{0,e_b-(q-u)\mathcal{E}\}.
 \end{eqnarray}
Furthermore, the average cost optimal policy $u^*$ satisfies
\begin{eqnarray}
\tilde{Z}(u^*)\le \beta\frac {\mathcal{E}}{\tau}p\le \tilde{Z}(u^*+1)
\end{eqnarray}
with
\begin{eqnarray}
\tilde{Z}(u)=\lim_{\alpha \to 1}Z(u).
\end{eqnarray}
\end{lemma}
\begin{IEEEproof}
See Appendix \ref{proof of Necessary condition for the discout optimality of u}.
\end{IEEEproof}
\begin{lemma}\label{discount optimal u for two special cases}
For $x=(q,a,e_b,e_a,p)$ satisfying
\begin{eqnarray}
Z\big(q-\min\{q,M\}\big)  > \beta \frac{\mathcal{E}}{\tau}p,
\end{eqnarray}
$u=q-\min\{q,M\}$ is the discount optimal policy.
In addition, for $(q,a,e_b,e_a,p)$ satisfying
\begin{eqnarray}
Z(q)< \beta \frac{\mathcal{E}}{\tau}p,
\end{eqnarray}
$u=q$ is the discount optimal policy.
\end{lemma}
\begin{IEEEproof}
See Appendix \ref{proof of discount optimal u for two special cases}.
\end{IEEEproof}
\emph{Remark: $u=q-\min\{q,M\}$, i.e., $k=\min\{q,M\}$ means charging as many EVs as possible. If the number of EVs in the queue is less than the charge point number $M$, charge all the EVs. Otherwise, charge $M$ EVs from the head of the queue. $u=q$, i.e., $k=0$ denotes charging no EV.}
\par
Based on Lemma \ref{discount optimal u for two special cases} and Lemma \ref{discout to UP}, we have
\begin{lemma}\label{average optimal u for two special cases}
For $x=(q,a,e_b,e_a,p)$ satisfying
\begin{eqnarray}
\tilde{Z}\big(q-\min\{q,M\}\big) > \beta \frac{\mathcal{E}}{\tau}p,
\end{eqnarray}
$u=q-\min\{q,M\}$ is the average cost optimal policy.
In addition, for $(q,a,e_b,e_a,p)$ satisfying
\begin{eqnarray}
\tilde{Z}(q)< \beta \frac{\mathcal{E}}{\tau}p,
\end{eqnarray}
$u=q$ is the average cost optimal policy.
\end{lemma}
\subsection{Two stationary deterministic policies}
Based on all above theoretical analysis, we propose the following two specific stationary deterministic policies.
For state $x=(q,a,e_b,e_a,p)$, we define the radical policy as $\Big(k=\min\{q,M\},w=\frac{\min\{e_b,k \mathcal{E}\}}{\tau}\Big).$ That is to say, we charge as many EVs as possible, and use the greedy policy for the battery energy allocation, i.e., if the required energy is not greater than the battery energy, then all the energy will be supplied
from the storage battery and no grid power will be used.
Otherwise, all the storage battery energy is allocated, and the
rest will be supplied from the power grid.
\par
In the radical policy, the average cost constraint is not considered. Then we propose another policy (i.e., the conservative policy) that guarantees the average cost constraint through satisfying the cost constraints in each period.
We call the policy
$\Big(k=\min\{q,M,\frac{e_b+\frac{\bar{B}}{p}\tau}{\mathcal{E}}\},w=\frac{\min\{e_b,k \mathcal{E}\}}{\tau}\Big)$ the conservative policy.
That is to say, we first guarantee that the cost of charging in each period is less than the average cost constraint, then charge as many EVs as possible and utilize the greedy policy for the battery energy allocation.
\par
In the whole paper, we assume that the power
from power grid and renewable energy generator is sufficient to stabilize the
queue length.
The stability issue such as the bounds
on average generation rate of renewable energy or average EV arrival rate
will be studied in future work.
\section{Numerical results}
In this section, we perform simulations to demonstrate the relations among the mean EV arrival, mean renewable energy arrival, upper bound of the average cost, average cost, and average EV queue length. Meanwhile, we consider different charge point numbers and capacities of the renewable energy storage battery.
In the simulations, the period length is $\tau=1$, and the size of the \lq\lq energy block\rq\rq ~is
$\mathcal{E}=10$.
\par
\begin{figure}[!t]
\centering
\subfigure[$M=50$]{\includegraphics[width=3.5in]{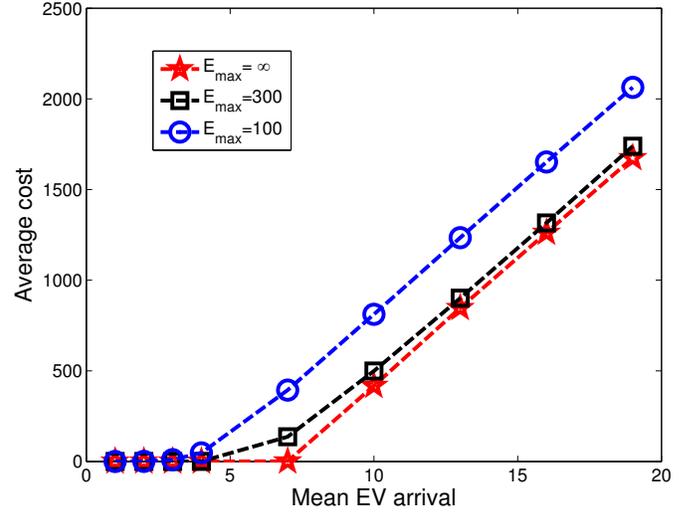}\label{mean EV arrivalM50}}
\subfigure[$M=8$]{\includegraphics[width=3.5in]{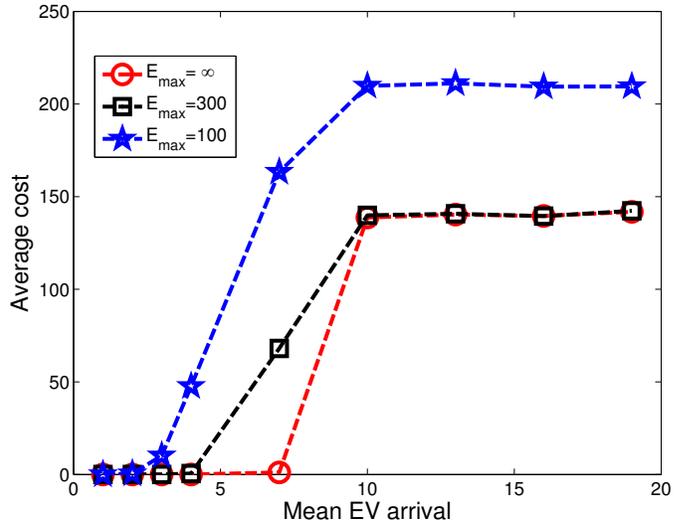}\label{mean EV arrivalM8}}
\caption{The average cost performance v.s. $\bar{A}$ under different values of $E_{max}$.}
\label{mean EV arrival}
\end{figure}

Fig. \ref{mean EV arrival} shows the average cost performance with respect to the mean EV arrival, $\bar{A}$. In the simulations, we utilize the radical policy.
We consider the i.i.d. cases of $A$, $E_a$, and $P$.
$A$ takes $0$ and $2 \bar{A}$ with equal probability $0.5$. $E_a$ takes values $\{0,50,100\}$ with probabilities $\{0.1,0.4,0.5\}$. $P$ takes values $\{5,10,20\}$ with probabilities $\{0.2,0.3,0.5\}$. The performance is averaged over $10^5$ periods. We set the number of charge points $M=50$ and $M=8$ in Fig. \ref{mean EV arrivalM50} and Fig. \ref{mean EV arrivalM8}, respectively.
Furthermore, we plot the curves for different storage battery capacities: $E_{max}=100$, $E_{max}=300$ and infinite capacity, respectively.
\par
In Fig. \ref{mean EV arrivalM50}, we can see that when $\bar{A}$ is small, the cost is nearly zero. However, when $\bar{A}$ is large (e.g., $\bar{A}\ge 10$), the cost increases rapidly with increase of $\bar{A}$ according to roughly a linear function.
It is because when $\bar{A}$ is small, the required energy is small and
the battery can supply the energy. Thus, no grid power will be consumed and the cost is zero. Once $\bar{A}$ is larger than a certain value, the required energy is larger than the battery energy, then the grid power will be utilized. As $M$ is large (compared to the considered $\bar{A}$), i.e., the restriction on the number of charge points will not influence the performance, we have $k=\min\{q,M\}=q$ with a high probability.
The grid power consumption will increase with increase of $\bar{A}$. Moreover, when $\bar{A}$ is large, the grid power becomes the main energy source. Based on (\ref{stage cost}), we derive that the cost varies with $\bar{A}$ roughly according to a linear relation.
\par
From Fig. \ref{mean EV arrivalM8}, we can find that the average cost is zero when $\bar{A}$ is small, and with increase of $\bar{A}$, the average cost increases. But once $\bar{A}$ is larger than a certain value, the average cost remains constant. It can be explained as follows: when $\bar{A}$ is small, the required energy can be supplied by the battery with a very high probability and no grid power is needed. Then the average cost is zero. When $\bar{A}$ increases, the required energy increases. Once the battery energy is not enough, the grid power will be consumed to fulfill the gap between the required energy and battery energy. With increase of $\bar{A}$, the grid power consumption increases since the average battery energy is constant. Thus, the average cost increases. However, when $\bar{A}$ is large enough, we get $k=\min\{q,M\}=M$ with a high probability because $M$ is not large in this simulations. Then, the required energy $k \times \mathcal{E}=M \times \mathcal{E}$, i.e., it becomes a constant. That means the grid power consumption is a constant also. Thus, the cost remains static.
\par
\begin{figure}[!t]
\centering
\includegraphics[width=3.5in]{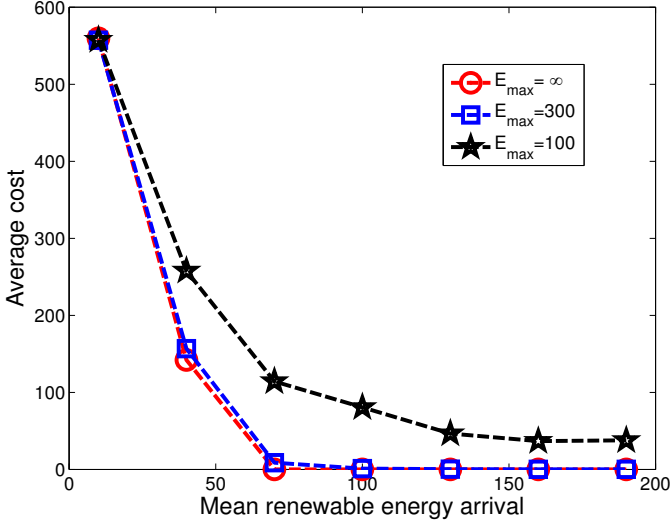}
\caption{The average cost under different $\bar{E_a}$ and $E_{max}$.}
\label{mean renewable energy arrival}
\end{figure}

Fig. \ref{mean renewable energy arrival} depicts the average cost performance with respect to the mean renewable energy arrival, $\bar{E_a}$. The radical policy is applied in the simulations. $A$ takes values $0$ and $10$ with equal probability $0.5$. $E_a$ take values $\{0,\frac{5}{7}\bar{E_a},\frac{10}{7}\bar{E_a}\}$ with probabilities $\{0.1,0.4,0.5\}$, respectively. $P$ is the same as in Fig. \ref{mean EV arrival} and $M=50$. $E_{max}=100$, $E_{max}=300$ and infinite capacity are also respectively considered in the simulations.
From the figure, we can find that the cost decreases with increase of $\bar{E_a}$. But once $\bar{E_a}$ is large enough, the cost almost remains static. First, in the range of small $\bar{E_a}$, when $\bar{E_a}$ increases, more free renewable energy will arrive and be stored in the battery. And then, the cost will decrease. If the battery capacity is large enough, all the arrived renewable energy can be stored in the battery. With the increase of  $\bar{E_a}$, the battery energy will increase all the time. Once the battery energy is larger than the required energy for charging, no grid power is needed then, and the cost becomes zero since that time. If the battery capacity is not large (e.g., $E_{max}=100$ in the figure), the overflow occurs when $\bar{E_a}$ is large. That is to say, the battery energy will remain $E_{max}$ even though we increase $\bar{E_a}$. On the other hand, $E_{max}$ is smaller than the required charge energy, so grid power is still needed. Consequently, the cost is non-zero and remains static.
\par
From Fig. \ref{mean EV arrival} and Fig. \ref{mean renewable energy arrival}, we can observe that the larger the battery capacity, the lower the cost. That is because when $E_{max}$ is larger, the probability of overflow will be lower (it is zero for infinite capacity). Then, less free renewable energy is wasted and the cost will be lower.
Furthermore, we can derive that if $\bar{A}$ is less a certain value or $\bar{E_a}$ is larger than a certain value, the average cost can be less than a certain value, Then we claim that when $\bar{A}$ is less a certain value or $\bar{E_a}$ is larger than a certain value, the radical policy is also optimal even when considering the constraint.\footnote{Notice that the radical policy is optimal for the mean EV queue delay minimization without the average cost constraint.}
\par
\begin{figure}[!t]
\centering
\includegraphics[width=3.5in]{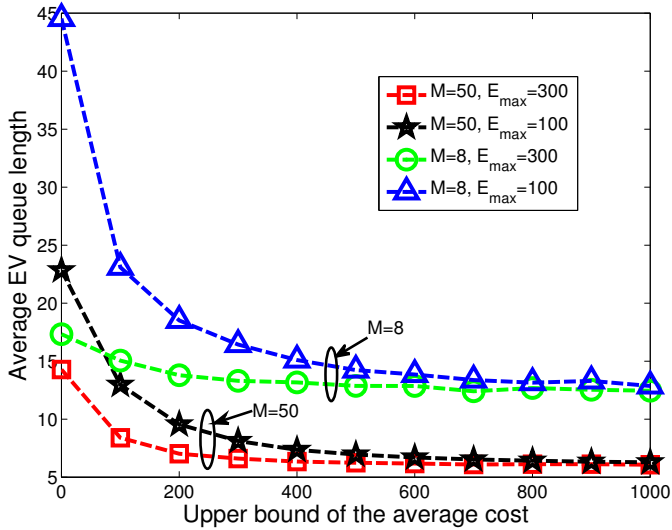}
\caption{The average EV length performance v.s. $\bar{B}$ .}
\label{Uppercostbound}
\end{figure}

Fig. \ref{Uppercostbound} illustrates the average EV queue length performance with respect to the upper bounds of the average cost when the conservative policy is applied. In the simulations, $A$ chooses values $\{0,12\}$ with equal probability $0.5$. $E_a$ and $P$ have the same settings as in Fig. \ref{mean EV arrival}. In the plotting, we consider different values of the battery capacity and charge point number. We can observe that the average length performance improves with increase of $\bar{B}$, and when $\bar{B}$ is larger than a certain value, the average length performance become almost constant. The reason is as follows: when $\bar{B}$ is small, $k=\min\{q,M,\frac{e_b+\frac{\bar{B}}{p}\tau}{\mathcal{E}}\} =\min\{q,\frac{e_b+\frac{\bar{B}}{p}\tau}{\mathcal{E}}\}$ with a high probability and it increases with increase of $\bar{B}$. Thus, the average EV queue length performance increases. Once $\bar{B}$ is large enough, we get $k=\min\{q,M\}$, and the average length remains static with respect to $\bar{B}$.
Additionally, by comparing the four curves, we can derive that the larger the capacity or the charge point number, the better the length performance.
\section{Conclusion}
We consider the scheduling of the EVs' charging at a charging station whose energy is provided from
both the power grid and local renewable energy. Under the uncertainty of the EV arrival, the renewable energy, the grid power price, and the charging energy of each EV, we study the mean delay optimal scheduling with the average cost constraint. We analyze the optimal policy
of the formulated MDP problem. In addition, two specific stationary policies (radical policy and conservative policy) are applied
in the simulations to reveal the impacts of relevant parameters on the performance.
%
%
\appendices
\section{Proof of Lemma \ref{equivalence energy demand queue and EV queue}}\label{proof of equivalence energy demand queue and EV queue}
First, the energy demand queue length and the EV queue length have the following relation.
$
Q_e[n]=\sum_{i=1}^{Q[n]}L_i
$
with $L_i$ being irrelevant to the queue state. Thus the average energy demand queue length is
$
\frac{1}{n}\sum_{j=1}^{n}Q_e[j]=\frac{1}{n}\sum_{j=1}^{n}\sum_{i=1}^{Q[j]}L_i^j.
$
Meanwhile, if an EV comes earlier than another EV, it will leave earlier in the EV queue serving. Using the queue mapping mechanism,
the earlier arrived EV will leave no later also in the energy demand queue serving.\footnote{Leave at the same time is possible.} That is to say, the queue mapping is an isotonic mapping.
Then, we claim that a policy minimizing the mean EV queue length results in minimal mean demand queue length, and vice versa.

\section{Proof of Lemma \ref{ConMDP to UnconMDP}}\label{proof of conMDP to unconMDP}
The proof is based on the results of \cite{CDC86:D. J. Ma A. M. Makowski and A. Shwartz}. We prove that for some $\beta$, the optimal policy $\pi^*$ of the unconstrained MDP (\ref{UP_beta}) (i.e., UP$_\beta$) satisfies 1) $\pi^*$ yields $B^{\pi^*}$ and $D^{\pi^*}$ as limits for all $x \in \mathcal{X}$;
2) $B^{\pi^*}=\bar{B}$. Observe that $\limsup$ and $\liminf$ are equal for each $\beta>0$ (since the controlled chain is ergodic and the policy is stationary \cite{Book02:E. Feinberg and A. Shwartz}).
\section{Proof of Lemma \ref{discout to UP}}\label{proof of discout to UP}
First, we derive that the conditions of Proposition 2.1 in \cite{MOR93:M. Schal} are satisfied. Then a discount optimal stationary policy exists.
Next, we prove that for some $x_0$, $V_{\alpha}(x)-V_{\alpha}(x_0)<\infty$. Third, there exits a policy $\pi \in \mathcal{A}$ and an initial state $x \in \mathcal{X}$ such that $J_\beta ^\pi <\infty$ in the practical problem. Otherwise, the cost is infinite for all policies and any policy is optimal. Accordingly, we can prove the lemma by applying Theorem 3.8 in \cite{MOR93:M. Schal}.
\section{Proof of Property \ref{Increasing in q}}\label{proof of Increasing in q}
We verify the increasing property by induction. According to (\ref{iteration for V_alpha}), $V_{\alpha,0}=0$ and $V_{\alpha,1}= \frac{\beta\big((q-\min\{q,M\})\mathcal{E}-e_b\big)^+p}{\tau}+ q$. The increasing property in $q$ holds. Assume $V_{\alpha,n-1}(q,a,e_b,e_a,p)$ is increasing in $q$. Depending on the values of $M$, we have the following two cases.
\par
Case 1: $M \ge q+1$.
Fix $(a,e_b,e_a,p)$, in the state $(q+1,a,e_b,e_a,p)$,
the set of feasible $u$ is $\{0,1,\cdots,q+1\}$ whereas it is $\{0,1,\cdots,q\}$ for state $(q,a,e_b,e_a,p)$.
Consider state $(q+1,a,e_b,e_a,p)$, let the optimal action be $(u^*,\eta^*)$ with $u^* \in \{0,1,\cdots,q\}$, hence
\begin{eqnarray}\label{1}
\lefteqn{
V_{\alpha,n}(q+1,a,e_b,e_a,p)  =
}
\nonumber\\
&&
\beta\Big(\frac{(q+1-u^*)\mathcal{E}}{\tau}-\frac{e_b-\eta^*}{\tau}\Big)^+p+(q+1)+\alpha\times
\nonumber\\
 && \mathbb{E}_{a,e_a,p}\left[V_{\alpha,n-1}(u^*+A,A,(\eta^*+E_a)^-,E_a,P)\right].
\end{eqnarray}
As $(u^*,\eta^*)$ is feasible in state $(q,a,e_b,e_a,p)$,
\begin{eqnarray}\label{2}
\lefteqn{
 V_{\alpha,n}(q,a,e_b,e_a,p) \le  \beta \Big(\frac{(q-u^*)\mathcal{E}}{\tau}-\frac{e_b-\eta^*}{\tau}\Big)^+p+ q
 }
  \nonumber\\
&+& \alpha \mathbb{E}_{a,e_a,p}\left[V_{\alpha,n-1}(u^*+A,A,(\eta^*+E_a)^-,E_a,P)\right]
 \nonumber\\
& \le& V_{\alpha,n}(q+1,a,e_b,e_a,p).
 \end{eqnarray}
If $(u^*,\eta^*)$ with $u^* =q+1$,
\begin{eqnarray}
\lefteqn{
V_{\alpha,n}(q+1,a,e_b,e_a,p)  =q +1
      }
\nonumber\\
&+&
\alpha \mathbb{E}_{a,e_a,p}\left[V_{\alpha,n-1}(q+1+A,A,(\eta^*+E_a)^-,E_a,P)\right].
\nonumber\\
\end{eqnarray}
Meanwhile, since $(q,\eta^*)$ is feasible in state $(q,a,e_b,e_a,p)$,
\begin{eqnarray}
   \lefteqn{
 V_{\alpha,n}(q,a,e_b,e_a,p) \le  q
 }
  \nonumber\\
 &+& \alpha \mathbb{E}_{a,e_a,p}\left[V_{\alpha,n-1}(q+A,A,(\eta^*+E_a)^-,E_a,P)\right]
 \nonumber\\
  &\stackrel{(a)}{\le} & V_{\alpha,n}(q+1,a,e_b,e_a,p),
 \end{eqnarray}
where (a) holds since the induction hypothesis.
 \par
 Case 2: $M \le q$.
The set of feasible $u$ is $\{0,1,\cdots,M\}$ in both the state $(q+1,a,e_b,e_a,p)$ and state $(q,a,e_b,e_a,p)$. Then we can prove the increasing property of $V_{\alpha,n}(q,a,e_b,e_a,p)$ by using (\ref{1}) and (\ref{2}).
\section{Proof of Property \ref{Non-increasing in eb}}\label{proof of Non-increasing in eb}
Based on (\ref{iteration for V_alpha}), the property can be proved through induction.
First, we have $V_{\alpha,0}=0$ and $V_{\alpha,1}= \frac{\beta\big((q-\min\{q,M\})\mathcal{E}-e_b\big)^+p}{\tau}+ q$. Thus the non-increasing property in $e_b$ holds for $n=0,1$. Next, assume $V_{\alpha,n-1}(q,a,e_b,e_a,p)$ is a non-increasing function of $e_b$. Fix $(q,a,e_a,p)$, for state $(q,a,e_b,e_a,p)$, let $(u^*,\eta^*)$ be the optimal policy. We get
\begin{eqnarray}\label{nonincreasing 1}
   \lefteqn{
 V_{\alpha,n}(q,a,e_b,e_a,p) =   \beta \Big(\frac{(q-u^*)\mathcal{E}}{\tau}-\frac{e_b-\eta^*}{\tau}\Big)^+p+
 }
 \nonumber\\
&&q+\alpha \mathbb{E}_{a,e_a,p}\left[V_{\alpha,n-1}(u^*+A,A,(\eta^*+E_a)^-,E_a,P)\right].  \nonumber\\
\end{eqnarray}
Since $(u^*,\eta^*)$ is feasible in state $(q,a,e_b+1,e_a,p)$, we derive
\begin{eqnarray}\label{nonincreasing 2}
   \lefteqn{
 V_{\alpha,n}(q,a,e_b+1,e_a,p)\le   \beta \Big(\frac{(q-u^*)\mathcal{E}}{\tau}-\frac{e_b+1-\eta^*}{\tau}\Big)^+p+
 }
 \nonumber\\
&&q+\alpha \mathbb{E}_{a,e_a,p}\left[V_{\alpha,n-1}(u^*+A,A,(\eta^*+E_a)^-,E_a,P)\right].  \nonumber\\
\end{eqnarray}
Combing (\ref{nonincreasing 1}) and (\ref{nonincreasing 2}), we get $$V_{\alpha,n}(q,a,e_b,e_a,p)\le  V_{\alpha,n}(q,a,e_b+1,e_a,p).$$
Then we complete the proof of the property.
\section{Proof of Property \ref{convexity}}\label{proof of convexity}
First, we prove the following proposition.
\begin{proposition}
For $\phi \in (0,1)$ and $\forall x_1,x_2,y$, we have $\phi\min\{x_1,y\}+(1-\phi)\min\{x_2,y\}\le \min\{\phi x_1+(1-\phi)x_2,y\}$.
\end{proposition}
\begin{IEEEproof}
The proposition can be verified by considering $\min\{x_1,x_2\}>y$, $\max\{x_1,x_2\}<y$, and $\min\{x_1,x_2\}\le y\le \max\{x_1,x_2\}$, respectively.
\end{IEEEproof}
The convexity is proved by induction.
For $n=0$, $V_{\alpha,0}=0$ and is convex. Assume $ V_{\alpha,n-1}(q,h,a,e_b,e)$ is convex in $(q,e_b)$. Fix $(q,a,e_b,e_a,p)$, let $(u_1,\eta_1)$ and $(u_2,\eta_2)$ be the optimal policy for $(q_1,e_{b1})$ and $(q_2,e_{b2})$. Then, we get
\begin{eqnarray}
  \lefteqn{
\phi V_{\alpha,n}( q_1,a, e_{b1},e_a,p)+(1-\phi)V_{\alpha,n}(q_2,a,e_{b2},e_a,p)
}
  \nonumber\\
&=& \phi\Big[\beta  \Big(\frac{(q_1-u_1)\mathcal{E}}{\tau}-\frac{e_{b1}-\eta_1}{\tau}\Big)p+q_1\Big]
\nonumber\\
&+&
(1-\phi)\Big[\beta\Big(\frac{(q_2-u_2)\mathcal{E}}{\tau}-\frac{e_{b2}-\eta_2}{\tau}\Big)p+ q_2\Big] \nonumber\\
&+& \alpha \mathbb{E}_{a,e_a,p}\Big[\phi V_{\alpha,n-1}(u_1+A,A,(\eta_1+E_a)^-,E_a,P)\nonumber\\
&+&(1-\phi)V_{\alpha,n-1}(u_2+A,A,(\eta_2+E_a)^-,E_a,P)\Big] \nonumber\\
&\stackrel{(b)}{\ge}& \beta\Big[\Big(\phi(q_1-u_1)+(1-\phi)(q_2-u_2)\Big)\mathcal{E}-\Big(\phi(e_{b1}-\eta_1)
\nonumber\\
&+&(1-\phi)(e_{b2}-\eta_2)\Big)\Big]\frac{p}{\tau}+[\phi q_1 \nonumber\\
&+&(1-\phi)q_2] +\alpha \mathbb{E}_{a,e_a,p}\Big[ V_{\alpha,n-1}(\phi u_1+(1-\phi)u_2\nonumber\\
&+&A,A,\phi (\eta_1+E_a)^-+(1-\phi)(\eta_2+E_a)^-,E_a,P)\Big]\nonumber\\
&\stackrel{(c)}{\ge}& \beta\Big[\Big(\phi(q_1-u_1)+(1-\phi)(q_2-u_2)\Big)\mathcal{E}-\Big(\phi(e_{b1}-\eta_1)
\nonumber\\
&+&(1-\phi)(e_{b2}-\eta_2)\Big)\Big]\frac{p}{\tau} \nonumber\\
&+&[\phi q_1+(1-\phi)q_2] +\alpha \mathbb{E}_{a,e_a,p}\Big[ V_{\alpha,n-1}(\phi u_1\nonumber\\
&+&(1-\phi)u_2+A,A,(\phi \eta_1+(1-\phi)\eta_2+E_a)^-,E_a,P)\Big]\nonumber\\
&\stackrel{(d)}{\ge}&
V_{\alpha,n}( \phi q_1+(1-\phi)q_2,a, \phi e_{b1}+(1-\phi) e_{b2},e_a,p), \nonumber
\end{eqnarray}
where (b) holds because of the convexity of $ V_{\alpha,n-1}(q,h,a,e_b,e)$, (c) holds because of Proposition 1 as well as Property \ref{Non-increasing in eb}, and
(d) holds since $(\phi u_1+(1-\phi)u_2,\phi \eta_1+(1-\phi)\eta_2)$ is feasible for $\phi( q_1,a, e_{b1},e_a,p)+(1-\phi)(q_2,a,e_{b2},e_a,p)$.
\section{Proof of Lemma \ref{Necessary condition for the discout optimality}} \label{proof of Necessary for the discout}
Let
\begin{eqnarray}
\lefteqn{
S(u,\eta)=\beta\Big(\frac{(q-u)\mathcal{E}}{\tau}-\frac{e_b-\eta}{\tau}\Big)p+ q
 }
 \nonumber\\
&+& \alpha \mathbb{E}_{a,e_a,p}\left[V_{\alpha}(u+A,A,(\eta+E_a)^-,E_a,P)\right].
\end{eqnarray}
First, we have
\begin{eqnarray}
 \lefteqn{
S(u+1,\eta)-S(u,\eta)=-\beta \frac{\mathcal{E}}{\tau}p
}
\nonumber\\
&+&\alpha \mathbb{E}_{a,e_a,p}\big[V_{\alpha}(u+1+A,A,(\eta+E_a)^-,E_a,P)
\nonumber\\
&-&V_{\alpha}(u+A,A,(\eta+E_a)^-,E_a,P)\big]
\end{eqnarray}
and
\begin{eqnarray}
  \lefteqn{
S(u-1,\eta)-S(u,\eta)=\beta \frac{\mathcal{E}}{\tau}p
}
\nonumber\\
&+&\alpha \mathbb{E}_{a,e_a,p}\big[V_{\alpha}(u-1+A,A,(\eta+E_a)^-,E_a,P)
\nonumber\\
&-&V_{\alpha}(u+A,A,(\eta+E_a)^-,E_a,P)\big].
\end{eqnarray}
Then applying $S(u^*+1,\eta^*)-S(u^*,\eta^*)\ge 0$ and $S(u^*-1,\eta^*)-S(u^*,\eta^*)\ge 0$, we obtain (\ref{u}).
Similarly, as
\begin{eqnarray}
 \lefteqn{
S(u,\eta+1)-S(u,\eta)=\beta \frac{p}{\tau}
}
\nonumber\\
&+&\alpha \mathbb{E}_{a,e_a,p}\big[V_{\alpha}(u+A,A,(\eta+1+E_a)^-,E_a,P)
\nonumber\\
&-&V_{\alpha}(u+A,A,(\eta+E_a)^-,E_a,P)\big]
\end{eqnarray}
and
\begin{eqnarray}
  \lefteqn{
S(u,\eta-1)-S(u,\eta)=\beta \frac{-p}{\tau}
}
\nonumber\\
&+&\alpha \mathbb{E}_{a,e_a,p}\big[V_{\alpha}(u+A,A,(\eta-1+E_a)^-,E_a,P)
\nonumber\\
&-&V_{\alpha}(u+A,A,(\eta+E_a)^-,E_a,P)\big],
\end{eqnarray}
we can reach (\ref{eta}) from $S(u^*,\eta^*+1)-S(u^*,\eta^*)\ge 0$ and $S(u^*,\eta^*-1)-S(u^*,\eta^*)\ge 0$.
In addition,
\begin{eqnarray}
  \lefteqn{
S(u+1,\eta+1)-S(u,\eta)=\beta \frac{p}{\tau}(1-\mathcal{E})
}
\nonumber\\
&+&\alpha\mathbb{E}_{a,e_a,p}\big[V_{\alpha}(u+1+A,A,(\eta+1+E_a)^-,E_a,P)
\nonumber\\
&-&V_{\alpha}(u+A,A,(\eta+E_a)^-,E_a,P)\big]
\end{eqnarray}
and
\begin{eqnarray}
 \lefteqn{
S(u-1,\eta-1)-S(u,\eta)=\beta \frac{p}{\tau}(\mathcal{E}-1)
}
\nonumber\\
&+&\alpha \mathbb{E}_{a,e_a,p}\big[V_{\alpha}(u-1+A,A,(\eta-1+E_a)^-,E_a,P)
\nonumber\\
&-&V_{\alpha}(u+A,A,(\eta+E_a)^-,E_a,P)\big].
\end{eqnarray}
Then, (\ref{u eta}) can be obtained by applying $S(u^*-1,\eta^*-1)-S(u^*,\eta^*)\ge 0$ and $S(u^*+1,\eta^*+1)-S(u^*,\eta^*)\ge 0$.
\section{Proof of Lemma \ref{Necessary condition for the discout optimality of u}} \label{proof of Necessary condition for the discout optimality of u}
First, based on Conjecture \ref{dimensionreduction}, we only need to consider the policy set
$\{(u,\eta):(u,\eta=\eta(u)) \cap (u,\eta) \ge (0,0)\} .$
Consequently,
\begin{eqnarray}
\lefteqn{
S(u,\eta(u))=\beta\Big(\frac{(q-u)\mathcal{E}}{\tau}-\frac{e_b-\eta(u)}{\tau}\Big)p+ q
}
 \nonumber\\
&+& \alpha \mathbb{E}_{a,e_a,p}\left[V_{\alpha}(u+A,A,(\eta(u)+E_a)^-,E_a,P)\right].
\nonumber\\
\end{eqnarray}
Then applying $S(u^*+1,\eta(u^*+1))-S(u^*,\eta(u^*))\ge 0$ and $S(u^*-1,\eta(u^*-1))-S(u^*,\eta(u^*))\ge 0$, we get
$
Z(u^*)\le \beta\frac {\mathcal{E}}{\tau}p\le Z(u^*+1).
$
Next, using Lemma \ref{discout to UP}, we reach the second half of the lemma.
\section{Proof of Lemma \ref{discount optimal u for two special cases}} \label{proof of discount optimal u for two special cases}
Following the proof of Lemma \ref{Necessary condition for the discout optimality of u},
we can prove the first half of the lemma
by contradiction.
Specifically, suppose $u=q-\min\{q,M\}$ is not the optimal solution, then $S(u^*-1,\eta(u^*-1))-S(u^*,\eta(u^*))\ge 0$ should hold. We have
$
Z\big(q-\min\{q,M\}\big) \le Z(u^*)
\le \beta \frac{\mathcal{E}}{\tau}p
$
and the contradiction occurs.
We can verify the second half of the lemma similarly by using contradiction. Assume $u=q$ is not the optimal solution, then $S(u^*+1,\eta(u^*+1))-S(u^*,\eta(u^*))\ge 0$ should be satisfied. Consequently, we get
$
Z(q) \ge
Z(u^*+1)
\ge \beta \frac{\mathcal{E}}{\tau}p.
$
The contradiction occurs then.



%

\end{document}